\documentclass{article}
\usepackage{wrapfig,amsmath,amssymb,amscd,theorem}
\usepackage{enumerate}
\usepackage{fancybox}
\usepackage[dvips]{graphicx,color}

\newcommand{\no}{\noindent}
\newtheorem{thm}{Theorem}[section]

\newtheorem{cor}[thm]{Corollary}

\newcommand{\R}{\mathbb{R}}

\newcommand{\pd}{\partial}

\renewcommand{\phi}{\varphi}
\renewcommand{\H}{\mathbb{H}}

\begin{document}

\title{\bf Extremal problems for the central projection}

\author{\bf Shigehiro Sakata}

\date{\today}

\maketitle

\def\c#1{\overset{\mbox{\tiny $\circ$}}{#1}}
\def\vect#1{\mbox{\boldmath $#1$}}
\def\vecttiny#1{\mbox{\tiny \boldmath $#1$}} 

\begin{abstract}
We consider a projection from the center of the unit sphere to a tangent space of it, the {\it central projection}, and study two area minimizing problems of the image of a closed subset in the sphere. One of the problems is the uniqueness of the tangent plane that minimizes the area for an arbitrary fixed subset. The other is the shape of the subset that minimizes the minimum value of the area. We also study the similar problems for the hyperbolic space.\\

\no{\it keywords and phrases}. central projection, extremal problem, area minimizing.
\end{abstract}

\section{Introduction}
Let $S^n$ denote the $n$-dimentional unit sphere, $S^n = \left\{ \left. x \in \R^{n+1} \right\vert \lvert x \rvert =1 \right\}$, $\sigma_n$ the spherical standard measure on $S^n$ and $\cdot$ the standard inner product of $\R^{n+1}$. Let $\Omega$ be an $n$-dimentional closed subset of $S^n$. The {\it polar set} of $\Omega$ is defined by
\[
\Omega^* = \bigcap_{y\in \Omega} \left\{\left. x \in S^n \right\vert x \cdot y \leq 0 \right\}.
\]
In the following, we always assume that $\Omega^*$ is not empty, i.e. $\Omega$ is contained in a hemisphere. Let $p_x$ ($x \in - \c{\Omega^*}$, where  $\c{\Omega^*}$ denotes the interior of the polar set $\Omega^*$) be the projection from $\Omega$ to the tangent space of $S^n$ at $x$ and $A_\Omega$ the map that assigns the area of $p_x (\Omega)$ to a point $x$ in $- \c{\Omega^*}$:
\[
p_{x} : \Omega \ni y \mapsto \frac{y}{x \cdot y} \in T_{x}S^n, \ \ A_\Omega (x) = {\rm Area}(p_x(\Omega)).
\] This projection $p_{x}$ is used for making a (local) world atlas and also called the {\it gnomonic projection} when $n=2$ in geography.
\begin{figure}[hbtp]
\centering
\scalebox{0.40}{\includegraphics[clip]{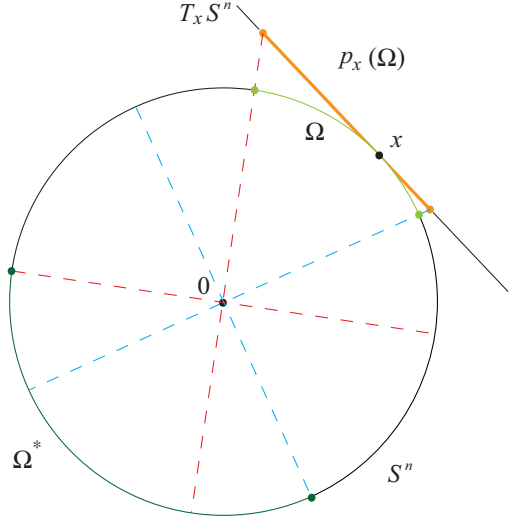}}
\caption{The polar set of $\Omega$ and $p_{x}(\Omega)$.}
\label{polar}
\end{figure}
Since the Jacobian of $p_{x}(y)$ is given by $\displaystyle Jp_{x}(y)=\frac{1}{(x \cdot y)^{n+1}}$, we have
\[
A_\Omega (x) = \int_{\Omega} \frac{1}{(x \cdot y)^{n+1}} d\sigma_n (y).
\]
\hspace{1.5em}F. Gao, D. Hug and R. Schneider showed the existence of a point that attains the minimum value of $A_\Omega$ and gave the characterization of it, but did not discuss the uniqueness in [1]. The existence of such a point follows from the fact that $A_\Omega$ is continuous on $- \c{\Omega^*}$ and diverges to $+ \infty$ as $x$ approaches the boundary of $- \Omega^*$. Direct calculation shows that if $x$ attains the minimum value of $A_\Omega$, then $x$ satisfies the following formula:
\[
\int_\Omega \frac{y}{(x \cdot y)^{n+2}} d\sigma_n (y) = A_\Omega (x) x.
\]\\  
\hspace{1.5em}In this paper we show the uniqueness of the $A_\Omega$ minimizer and estimate the minimum value of $A_\Omega$ for a closed subset $\Omega$ in $S^n$ having the same area with a disc in $S^n$. We also study the similar problems for the hyperbolic space.\\
 
\no{\bf Acknowledment.}
The author would like to express his deep gratitude to his advisors Jun O'Hara and Gil Solanes for giving many kind advices to him.
\section{Extremal value of {\boldmath $A_\Omega(x)$}}
Let $\Omega$ and $A_\Omega$ be as in the introduction. 
\begin{thm}\label{main1} $A_\Omega$ has a unique minimum point. 
\end{thm}
{\bf Proof.} By using the polar coordinate, we can put 
\begin{align*}
x &= x(\theta_1, \ldots, \theta_n)\\
  &= \left( \cos \theta_1, \sin \theta_1 \cos \theta_2, \ldots, \sin \theta_1 \cdots \sin \theta_{n-1} \cos \theta_n, \sin \theta_1 \cdots \sin \theta_{n-1} \sin \theta_n   \right),  
\end{align*} where $(\theta_1, \ldots, \theta_{n-1}, \theta_n) \in [0, \pi ] \times \cdots \times [0, \pi ] \times [0, 2\pi ]$. Since $\displaystyle \frac{\pd^2 x}{\pd \theta_1^2} = -x$, we have
\[
\frac{\pd^2 A_\Omega}{\pd \theta_1^2} (x) = (n+1) \int_{\Omega} \frac{(x \cdot y)^2 + (n+2) (\frac{\pd x}{\pd \theta_1} \cdot y)^2}{(x \cdot y)^{n+3}} d\sigma_n (y) >0
\] on $- \c{\Omega^*}$. Suppose that there exist two points $x'$ and $x''$ that attain the minimum value of $A_\Omega$. By a rotation of $S^n$, we may assume that the sub-arc between $x'$ and $x''$ of the great circle constant $\theta_i$ coordinates for $i \neq 1$. Then we have
\[
\frac{\pd A_\Omega}{\pd \theta_1} (x') =\frac{\pd A_\Omega}{\pd \theta_1} (x'') =0,
\] which is a contradiction.\hspace{\fill}$\Box$
\begin{cor}\label{cor} If $\Omega$ is point symmetric at $x_s$, then $x_s$ is the unique minimum point of $A_\Omega$.
\end{cor}

By Corollary \ref{cor}, the center of a disc $D$ in $S^n$ with non-empty polar set is the unique minimum point of $A_D$. This fact can also be indicated by the {\it moving plane method} (cf. [2]) became of the symmetry of the integrand of the partial derivative of $A_\Omega$.
\begin{thm}\label{main2} Let $D$ be a disc in $S^n$ with non-empty polar set. If $\sigma_n(\Omega) = \sigma_n(D)$, then 
\[
\min_{x \in - \c{D^*}} A_D(x) \leq \min_{x \in - \c{\Omega^*}} A_\Omega (x)
\]and that equality holds if and only if $\Omega$ is a spherical cap.
\end{thm}
{\bf Proof.} By a translation of $S^n$, we may assume that the center of $D$ coincides with the (unique) minimum point of $A_\Omega$. Let $x_c$ denote the center of $D$. Then we obtain
\[
A_\Omega (x_c) - A_D(x_c) =\int_{\Omega \backslash (\Omega \cap D)} \frac{1}{(x_c \cdot y)^{n+1}} d\sigma_n(y)- \int_{D \backslash (\Omega \cap D)} \frac{1}{(x_c \cdot y)^{n+1}} d\sigma_n(y) \geq 0
\] since the following inequality holds for $y' \in \Omega \backslash (\Omega \cap D)$ and $y'' \in D \backslash (\Omega \cap D)$:
\[
x_c \cdot y' = \cos \angle (x_c, y') \leq \cos \angle (x_c, y'') = x_c \cdot y'.
\] That equality holds if and only if $\sigma_n(\Omega \backslash (\Omega \cap D))=0$, namely, $\Omega$ is a disc in $S^n$. \hspace{\fill}$\Box$\\

We can see the following theorem with the same argument of Theorem \ref{main2}.
\begin{thm}\label{main2'}
Let $D$ be a disc in $S^n$ and $f:[0, +\infty) \to \R$ a strictly increasing continuous function. For an $n$-dimentional closed subset $K$ in $S^n$, if $\sigma_n(K) = \sigma_n(D)$, then we have
\[
\min_{x \in S^n}\int_D f({\rm dist}_{S^n}(x,y))d \sigma_n (y) \leq \min_{x \in S^n}\int_K f({\rm dist}_{S^n}(x,y))d \sigma_n (y).
\]
\end{thm}
\section{Hyperbolic case}
Let $\langle \cdot , \cdot \rangle$ denote the indefinite inner product of $\R^{n+1}$ given by $\langle x, y \rangle = x_1y_1 + \cdots + x_ny_n - x_{n+1}y_{n+1}$, $\R_1^{n+1}$ the $(n+1)$-dimensional Euclidean space with $\langle \cdot, \cdot \rangle$, and $\mathbb{H}^n$ the Lorents model of the $n$-dimensional hyperbolic space, $\mathbb{H}^n = \left\{ \left. x \in \R_1^{n+1} \right\vert \langle x , x \rangle = -1, x_{n+1} > 0 \right\}$. Let $\Omega$ be an $n$-dimensional compact subset in $\mathbb{H}^n$ with non-empty interior. We define the maps $p_x$ ($x \in \mathbb{H}^n$) and $A_\Omega$ by the similar way with the spherical case:
\[
p_x: \Omega \ni y \mapsto -\frac{y}{\langle x , y \rangle} \in T_x\mathbb{H}^n, \ \ A_\Omega(x) = {\rm Area}(p_x(\Omega )).
\] Since the Jacobian of $p_x(y)$ is given by $\displaystyle Jp_x(y) = \frac{1}{(- \langle x , y \rangle)^{n+1}}$, we have
\[
A_\Omega (x) = \int_\Omega \frac{1}{(- \langle x,y \rangle )^{n+1}} d\mu_n (y),
\] where $\mu_n$ is the standard hyperbolic measure on $\mathbb{H}^n$.\\
\hspace{1.5em}The existence of a $A_\Omega$ maximizer follows from the fact that $A_\Omega$ is continuous on $\mathbb{H}^n$ and converges to 0 as $x_{n+1}$ goes to $+\infty$. The simiraly computation with [1] shows that if $x$ attains the maximum value of $A_\Omega$, then $x$ satisfies the following formula:
\[
\int_\Omega \frac{y}{(- \langle x, y \rangle )^{n+2}} d\mu_n (y) = A_\Omega (x) x.
\] The uniqueness of such a point does not always hold. For example, the following set $\Omega$ in $\mathbb{H}^1$ has two maximum points of $A_\Omega$:
\[
\Omega = \left\{\left. (\sinh \theta,\cosh \theta) \right\vert -2 \leq \theta \leq -1 \right\} \cup \left\{\left. (\sinh \theta, \cosh \theta) \right\vert 1 \leq \theta \leq 2 \right\}.
\]
\hspace{1.5em}On the other hand, we can obtain the following theorems corresponding to Theorem \ref{main2} and \ref{main2'} with the same arguments.
\begin{thm}\label{main3} Let $D$ be a disc in $\mathbb{H}^n$. If $\mu_n(\Omega) = \mu_n(D)$, then we have
\[
\max_{x\in \mathbb{H}^n} A_\Omega (x) \leq \max_{x\in \mathbb{H}^n} A_D (x)
\] and that equality holds if and only if $\Omega$ is a disc in $\mathbb{H}^n$.
\end{thm}
\begin{thm}\label{main3'}
Let $D$ be a disc in $\H^n$ and $f:[0, +\infty) \to \R$ a strictly decreasing continuous function. For an $n$-dimentional compact subset $K$ in $\H^n$, if $\mu_n(K) = \mu_n(D)$, then we have
\[
\max_{x \in \H^n}\int_K f({\rm dist}_{\H^n}(x,y))d \mu_n (y) \leq \max_{x \in \H^n}\int_D f({\rm dist}_{\H^n}(x,y))d \mu_n (y).
\]
\end{thm}

\no 
Department of Mathematics and Information Science,\\
Tokyo Metropolitan University,\\
1-1 Minami Osawa, Hachiouji-Shi, Tokyo 192-0397, Japan\\
E-mail: sakata-shigehiro@ed.tmu.ac.jp

\end{document}